\documentclass[a4paper,fleqn]{cas-sc}

\usepackage[numbers, sort&compress]{natbib}
\usepackage{amssymb}
\usepackage{amsmath}
\usepackage{float}
\usepackage{hyperref}

\usepackage{graphicx}
\usepackage{caption}
\usepackage{subcaption}

\usepackage{xcolor}
\usepackage{textcmds}
\usepackage{svg}
\usepackage{lipsum}

\def\tsc#1{\csdef{#1}{\textsc{\lowercase{#1}}\xspace}}
\tsc{WGM}
\tsc{QE}

\begin{document}
\let\WriteBookmarks\relax
\def\floatpagepagefraction{1}
\def\textpagefraction{.001}

\definecolor{blue}{rgb}{0,0,1}
\newcommand{\IEE}[1]{\textcolor{blue}{#1}}
\definecolor{red}{rgb}{1,0,0}
\newcommand{\REV}[1]{\textcolor{red}{#1}}

\shorttitle{}    

\shortauthors{}  

\title[mode = title]{Implementing Dynamic Power Feed-In Limitations of Photovoltaic Systems in Distribution Grids for Generation Expansion Planning}





%

\author[1,2]{Alexander Konrad}[type=editor,
        orcid=0000-0001-9355-0746]

\cormark[1]


\ead{alexander.konrad@tugraz.at}

\ead[url]{www.tugraz.at/institute/iee/home}

\credit{Conceptualization of this study, Methodology, Software, Validation, Writing - Original draft preparation}

\affiliation[1]{organization={Institute of Electricity Economics and Energy Innovation, Graz University of Technology},
            addressline={Inffeldgasse 18},
            city={Graz},
            postcode={8010}, 
            state={Styria},
            country={Austria}}

\affiliation[2]{organization={Research Center ENERGETIC, Graz University of Technology},
            addressline={Rechbauerstraße 12}, 
            city={Graz},
            postcode={8010}, 
            state={Styria},
            country={Austria}}            

\author[1,2]{Robert Gaugl}[type=editor,
        orcid=0000-0003-4112-4483]
\credit{Conceptualization, Software, Review}
\author[3]{Christoph Maier}[type=editor]
\credit{Data curation, Validation}
\author[1,2]{Sonja Wogrin}[type=editor,
        orcid=0000-0002-3889-7197] 
\credit{Supervision, Software, Review}
        




\affiliation[3]{organization={Netz Niederösterreich GmbH},
            addressline={EVN Platz 1}, 
            city={Maria Enzersdorf},
            postcode={2344}, 
            state={Lower Austria},
            country={Austria}}

\cortext[1]{Corresponding author}




\begin{abstract}
The rapid growth of photovoltaic (PV) systems in Austria's medium- and low-voltage grids has intensified challenges in grid access, with technical limits increasingly leading to restrictions on full feed-in power. This issue has sparked discussions about limiting PV feed-in power and the implications for both generated and curtailed PV energy. At the same time, expanding PV capacity remains critical to achieving future climate targets. However, there is a lack of robust methodologies of quantify the impact of PV feed-in limitations when implemented in an optimization model. This impact affects both the curtailed energy and the increase in maximum PV installation capacity and total energy production. To address this gap, we have developed a mathematical formulation of dynamic PV feed-in limitations and integrated it into an optimization model. This approach enables a comprehensive evaluation of its effects on PV integration potential and energy curtailment, validated through case studies on four representative real-world Austrian medium- and low-voltage grids. We analyzed maximum PV expansion, energy generation, and curtailment under feed-in constraints. The results highlight the potential for integrating up to 32\% additional PV systems within existing infrastructure while keeping PV curtailment relatively low, i.e. at 2\%. We provide actionable insights for grid operators and policymakers aiming to balance renewable energy expansion with grid reliability.
\end{abstract}

\begin{graphicalabstract}
    \begin{figure}[h!]
        \centering
        \includegraphics[width=1\textwidth]{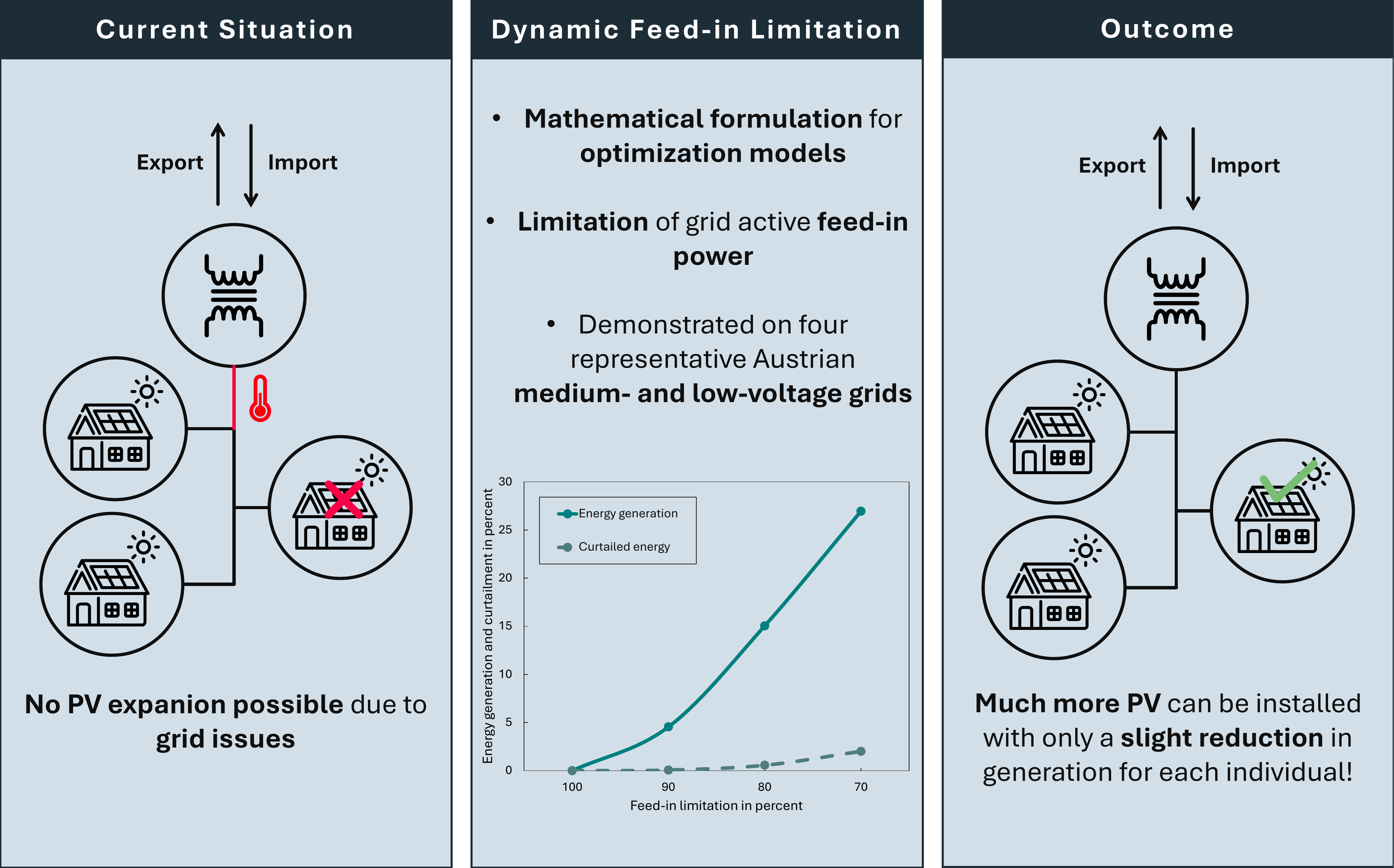}
    \end{figure}
\end{graphicalabstract}

\begin{highlights}
\item By means of dynamic feed-in limitation PV capacity can be increased by up to 32\%.
\item Results indicate feed-in limitations cause an overall energy loss of only 2\%.
\item An approach to modeling and analyzing dynamic feed-in limitation for PV expansion.
\item Methodology demonstrated on characteristic Austrian low- and medium-voltage grids.
\end{highlights}

\begin{keywords}
distribution systems\sep active-power curtailment\sep solar power\sep energy system optimization\sep PV feed-in limitation\sep
\end{keywords}

\maketitle

\begin{table}[h!]
\begin{tabular}{l l}
    \textbf{Nomenclature:}\\
    Sets:\\
    $h$                     & Time periods (usually hours)\\
    $i$                     & Bus of distribution grid\\
    $g$                     & Generation units\\
    $pv(g)$                 & Subset of photovoltaic renewable generation units\\
    $gi(g,i)$               & Generator g connected to node i\\
    \\
    Parameters:\\
    $FL$                    & Feed-in power limitation; maximum feed-in depending on installed power (\%)\\
    $D_{h,i}^P$             & Active power demand (MW)\\
    $\overline{P}_g$        & Technical maximum power of unit (MW)\\
    $CF_{h,i,g}$            & Renewable profile per technology and node; Capacity factor (p.u.)\\
    $C^{Imp}$               & Cost of imports (€/MWh)\\
    $C^{Exp}$               & Cost of exports (€/MWh)\\
    $C^{PNS/QNS}$           & Cost of power non-served (€/MWh)/(€/MVArh)\\
    $C^{EPS/EQS}$           & Cost of excess power served (€/MWh)/(€/MVArh)\\
    $M$                     & Big positive constant for activating or deactivating constraints\\
    \\
    Variables:\\
    $p_{h,g}$               & Active power generation of the unit (MW)\\
    $sp_{h,g}$              & Curtailment of renewable generation units (MW)\\
    $p_{h}^{Imp/Exp}$       & Active power import/export (MW)\\
    $q_{h}^{Imp/Exp}$       & Reactive power import/export (MVAr)\\
    $\alpha_{h,i}$          & Allows curtailment \{0,1\}\\
    $scal$                  & Scaling variable for all scalable PV systems (positive continuous)\\
    $pns_{h,i} / qns_{h,i}$ & Active / reactive power non-served (MW / MVAr)\\
    $eps_{h,i} / eqs_{h,i}$ & Excess active / reactive power served (MW / MVAr)\\
\end{tabular}
\end{table}

\section{Introduction}
\label{sec:Introduction}

The transition of the electricity sector from fossil fuels to renewable energy sources and towards complete decarbonization is one of the biggest challenges of our time. The European Green Deal targets a reduction in greenhouse gas emissions of at least 55\% by 2030, 90\% by 2040 compared with 1990 levels, and net zero by 2050 \cite{EC2020}. Austria has set even more ambitious targets and aims to achieve climate neutrality by 2040~\cite{OesterreichGV2024}. To achieve this goal, the Renewable Energy Expansion Act (EAG) -- a law in force since 2021 -- aims to generate an additional 27~TWh from renewable sources by 2030 compared with 2020 \cite{EAG2021}. Accordingly, the 55.4~TWh already generated from renewable sources in 2020~\cite{StatistikAustriaEnergiebilanz}, this corresponds to covering Austria's annual electricity demand net-nationally, i.e. on an annual average. Of this 27~TWh target, 11~TWh will come from photovoltaic (PV) alone, which corresponds to a total installed capacity of 13~GW of PV, including the 2~GW already installed in 2020~\cite{EAG2021}. This installed PV capacity should increase to 41~GW by 2040~\cite{OENIP2024}. PV is one of the most important and promising technologies for the United Nations' Sustainable Development Goal of affordable and clean energy~\cite{UN_SDG}. PV is already one of the energy sources with the lowest generation costs~\cite{FrauenhoferISE2024}. Furthermore, it is possible to install small systems on almost any building and thus generate extremely cheap and clean energy. The installation of PV systems can be executed in a variety of ways, including conventional mounting on rooftops or ground-mounting larger systems. However, it should be noted that there are also exceptional installation variants, such as facade integration, roof-integrated systems, and Solar-Photovoltaic-Trees~\cite{Gangwar2019,Kumar2019}.

Most PV systems in Austria are installed in a decentralized manner on rooftops, whereas somewhat larger systems are predominantly ground-mounted PV systems~\cite{E-Control2024}. However, in contrast to most other renewable energy sources, such as hydropower and wind power, these PV systems are predominantly (>99\% of all systems) located at grid levels 5 to 7 (medium- and low-voltage levels)~\cite{E-Control2024}. PV expansion is creating more and more congestions at these grid levels because they were originally built as \qq{one-way street,} to cover consumption, not to receive PV feed-in. These bottlenecks, which limit the further expansion of PV systems, manifest themselves at both medium- and low-voltage levels in the form of voltage problems or by exceeding the thermal limit of the infrastructure~\cite{Aziz2017,Gandhi2020}.

To integrate a new PV system into the Austrian distribution grid, an inquiry must first be submitted to the local distribution system operator (DSO). In most cases, only a \qq{static grid calculation} is carried out by the DSO, with the lowest grid load with simultaneous maximum PV generation being assumed (low-load case) and checked for compliance with all limit values (voltage and thermal power limits) using an AC power flow~\cite{OE2020}. A few DSOs are now working with the real time series for demand and generation from corresponding grid~\cite{OE2020, E-Control2023}.
Flexible grid access for feeders is currently being discussed in Austria in a draft of the new Electricity Economic Act (ELWG~\cite{ELWG2024}). This provides the possibility for DSOs to impose a static limitation on generation capacity to a maximum of 80\%. This is to be achieved without consideration of self-consumption, thereby also establishing a dynamic limitation on generation. However, the DSO must fix the bottleneck and the associated limitation of maximum feed-in power for a duration of 6 to 18 months (depending on the grid level).
The approach with static and dynamic curtailment bears a resemblance to the feed-in limitation in Germany, which was recently overturned as part of the amendment to the Renewable Energy Sources Act (EEG amendment). In contrast to the Austrian approach, the limit for PV systems up to 25 kW was fixed at 70\% and was not subject to temporal limitations. Another jurisdiction that employs a feed-in limit is the Australian state of South Australia, where the feed-in power for new systems has been reduced to 1.5~kW. However, there is the possibility of feeding up to 10~kW into the distribution grid of South Australian Power Networks (SAPN). PV system operators are obligated to use inverters that are capable of receiving the maximum feed-in power applicable to them from the DSO on a quarterly basis, depending on the status of the grid.

A review of the literature reveals a multitude of approaches that address the challenges of integrating PV into medium- and low-voltage grids, as well as potential solutions. In most cases, compliance with voltage limits has been identified as the most significant challenge~\cite{Aziz2017,Gandhi2020,Haghighian2021,Malekpour2017,Mukherjee2024,Ozdemir2023,Stetz2013,Sun2022}. Numerous studies have developed methodologies and algorithms to address the issue of voltage violations when large amounts of PV are added to the grid. A synthesis of the extant literature shows that there are four main approaches. The most obvious approach to reduce congestion is traditional grid expansion (i). To determine the optimal grid expansion for distribution grids, researchers developed a planning method in~\cite{Laribi2021}. In addition to the reinforcement of power lines and transformers, the use of adjustable local grid transformers (OLTCs) has emerged as an effective method, particularly for addressing voltage issues~\cite{Aziz2017,Sun2022,Stetz2013}. Option (ii) controls reactive power and offers various control options. By controlling reactive power in a targeted manner, it is possible to reduce over-voltages without reducing apparent power generation, thus remaining within voltage limits~\cite{Aziz2017,Haghighian2021,Malekpour2017}. Limiting active power (iii) also results in a reduction of over-voltages, in contrast to reactive power control, energy is curtailed.~\cite{Aziz2017,Laribi2021,Malekpour2017,Liu2020,Ozdemir2023,Stetz2013,Luthander2017}. Another option (iv) is the use of battery energy storage systems (BESS)~\cite{Laribi2021,Aziz2017,Moshoevel2014,Petrou2021,Luthander2017}. However, it is important that BESS are charged intelligently, as otherwise, the BESS may already be at maximum state of charge during periods of maximum generation, negating any potential added value.
In the optimization distribution grid planning method in~\cite{Laribi2021}, a cost comparison was made between grid expansion, BESS investments, and dynamic curtailment of active power. The results indicated that the curtailment of active power was the most cost-effective option. In the study by Malekput and Pahwa~\cite{Malekpour2017}, they compared different variants of reactive power control (e.g., $\cos\varphi(P)$-, $Q(V)$-, and $Q(V)/P(V)$-control) with active power limitation. Pure active power limitation demonstrates the lowest active power losses following the $Q(V)$-control. Unlike most of the other referenced studies, Stetz et al.~\cite{Stetz2013} increased the installed PV power in various simulations. Other studies generally assumed an existing bottleneck, which was effectively removed using different methods.

To our knowledge, there is no existing mathematical formulation of this PV feed-in limitation in the literature. Unlike previous studies, our model focuses specifically on maximizing the potential installed capacity of PV systems. Moreover, while some of the cited studies have also considered fairness across different PV systems~\cite{Liu2020, Mukherjee2024, Sun2022}, our approach addresses this uniquely by implementing the dynamic feed-in limitation uniformly across all PV systems, regardless of their grid location. The original contribution of this paper lies in developing an optimization model grounded in the concept of \qq{dynamic feed-in limitation,} a topic currently under discussion in Austria. In order to quantify the impact of this method on the potential increase in PV system capacity, energy generation, and curtailed energy, the model was applied to four real-world grids for validation. This analysis produced case studies that offer insights into the practical feasibility and benefits of the proposed approach.

The remainder of this paper is structured as follows: Section~\ref{sec:Methodology} contains a detailed description of the methodology used, as well as a description of the optimization model (digital representation of the medium-voltage and low-voltage electricity system) used to carry out this work. Section~\ref{sec:Case-Study} deals with the case studies and their data. Section~\ref{Sec:Results} describes the results of the case studies, which include curtailed energy, the additional installed PV capacities, the energy from additional generation, and bottlenecks identified in each grid. Finally, section~\ref{sec:Conclusion} summarizes and concludes this paper.

\section{Methodology}
\label{sec:Methodology}
This section describes the methodology employed to assess the dynamic feed-in limitation of PV systems in medium- and low-voltage grids. First, in subsection~\ref{subsec:General-methodology} we provide a detailed description of the dynamic feed-in limitation discussed in Austria. Then, in subsection~\ref{subsec:Model-and-math-formulation} we present a mathematical formulation of this dynamic feed-in limitation that can be included in power system optimization models, e.g. generation and transmission expansion planning models.

\subsection{General methodology}
\label{subsec:General-methodology}

\begin{figure}[h!]
\centering
\includegraphics[width=84 mm]{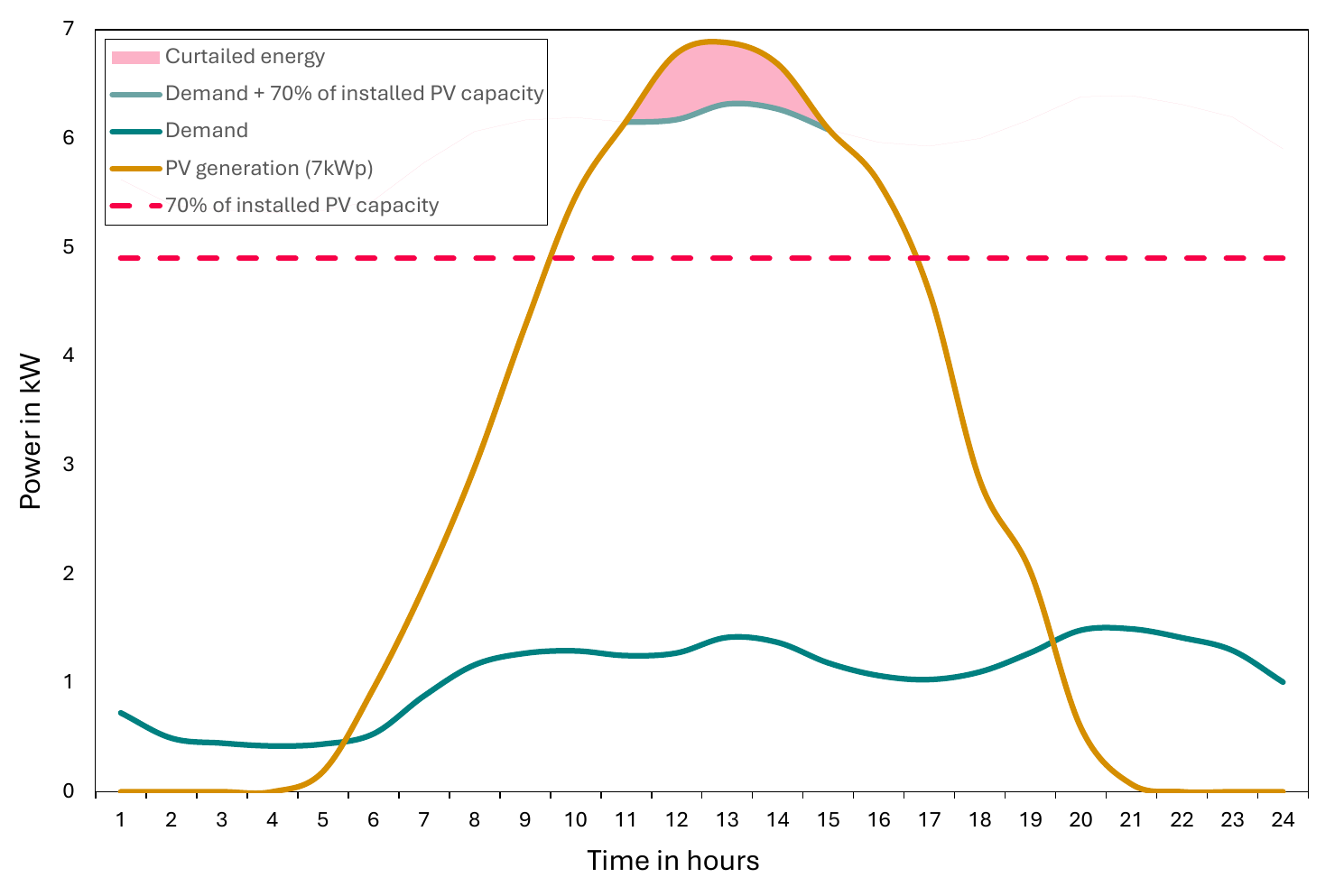}
\caption{Stylized example of dynamic feed-in limitation of PV.}\label{fig:Dyn-feed-in-limitation}
\end{figure}

The general idea of dynamic feed-in limitation is to reduce grid-active feed-in power to the medium- or low-voltage grid rather than limiting a PV system's generation capacity. This can be achieved by taking into account not only PV generation but also concurrent self-consumption. The following provides a small numerical example to showcase the dynamic feed-in limitation. First, a regulatory entity sets a percentage of the installed maximum peak power, e.g. 70\%, representing the maximum feed-in of a PV system. For a 7~kWp PV system, the maximum grid-active feed-in power would be 4.9~kW (70\% of 7~kW module peak power). If the limitation is static, 2.1~kW would be lost at maximum generation. With dynamic feed-in limitation, however, the self-consumption can influence this. If self-consumption is 1.4~kW at this moment, then only 0.7~kW would be curtailed. Figure~\ref{fig:Dyn-feed-in-limitation} illustrates this simple example for a high-yield day. The dashed line represents the feed-in limitation of 70\% of the installed PV capacity. If self-consumption or demand (dark-green line) were to be ignored, then all PV generation above the dashed line would be curtailed (static limitation). However, the dynamic feed-in limitation accounts for the demand and subtracts it from the PV generation. Hence, there are only a few hours (light red area) where PV generation minus demand exceeds the 70\% limitation and is therefore curtailed.

The method employed in this paper does not consider storage systems. The strategic implementation of BESS enables the temporary storage of curtailed energy for later use. However, in scenarios where storage systems are not utilized judiciously, it is conceivable that these systems have already been fully charged, rendering them incapable of storing additional energy during periods of curtailment. The feasibility of achieving near-zero curtailed energy output depends on factors including PV generation, battery storage capacity, and the system's operational mode.

\subsection{Optimization model and mathematical formulation}
\label{subsec:Model-and-math-formulation}

In this section, we propose a mathematical formulation of the dynamic feed-in limitation for PV systems, which can be included in standard power system optimization models. For the case studies in this paper, we have used the Low-carbon Expansion Generation Optimization (LEGO) model~\cite{Wogrin2022} to which we have added the formulation of the dynamic feed-in limitation. The original model is open-source and available on GitHub\footnote{\url{https://github.com/IEE-TUGraz/LEGO}}.
For the sake of completeness, in the following we provide a brief description of the LEGO model. LEGO is a highly modular, flexible and deterministic mixed-integer, quadratically constrained optimization model (cost minimization) for investigating a variety of aspects in the energy sector. These comprise power plant dispatch decisions, including ramp-up and ramp-down for thermal power plants, load flow calculations based on a DC-optimal power flow (DC-OPF) or a second-order cone approximation of the AC-optimal power flow (AC-OPF), long-term expansion planning of generation and transmission infrastructure, etc. For a detailed description of LEGO and its mathematical formulation, interested readers are referred to~\cite{Wogrin2020},\cite{Wogrin2022}.\\

We now describe the formulation for the dynamic feed-in limitation of PV, which represents the original contribution of this paper. The objective function~\eqref{eq:ObjFunc} aims to minimize total system costs, i.e. costs incurred for the exchange of active and reactive power $p_{h}^{Imp/Exp}$, $ q_{h}^{Imp/Exp}$ with the higher-level grid (medium-voltage or high-voltage grid). In the event that consumption cannot be covered or a given generation cannot be transported away from a node, this energy can be covered or removed via the variables $pns_{h,i} / qns_{h,i}$ and $eps_{h,i} / eqs_{h,i}$. However, this approach is associated with exorbitantly high costs and is considered by the model as a last resort. The model parameters are shown in Table~\ref{tbl:Parameter}.

\begin{equation}\label{eq:ObjFunc}
    \begin{aligned}
        \min \sum_{h} (p_{h}^{\text{Imp}} + q_{h}^{\text{Imp}}) \cdot h \cdot C^{\text{Imp}} 
        &- \sum_{h} (p_{h}^{\text{Exp}} + q_{h}^{\text{Exp}}) \cdot h \cdot C^{\text{Exp}} \\
        &+ \sum_{h} \sum_{i} (pns_{h,i} + qns_{h,i}) \cdot h \cdot C^{\text{PNS/QNS}} \\
        &+ \sum_{h} \sum_{i} (eps_{h,i} + eqs_{h,i}) \cdot h \cdot C^{\text{EPS/EQS}}
    \end{aligned}
\end{equation}

Constraints~\eqref{eq:MaxProd}-\eqref{eq:Regulation3} describe the dynamic feed-in limitation for PV systems. Note that~\eqref{eq:MaxProd}-\eqref{eq:Regulation3} only determine active power. This is because the DSO specified $\cos\varphi = 1$ for PV generation and consumption. However, the formulation can easily be extended to also include reactive power. Constraint~\eqref{eq:MaxProd} defines the upper and lower bounds of PV generation $p_{h,pv}$. Constraint~\eqref{eq:Gen+Curt} ensures that generated and curtailed energy $sp_{h,pv}$ of each PV generator correspond to the maximum possible generation in each time step $h$. The positive variable $scal$ is used to scale the generation of new PV systems up to the limits of the grid. This increases the PV generation and potential exports, which in turn reduces the objective function. Note that in order to not disadvantage individual installations, $scal$ represents the system-wide scaling of PV capacity. Thus, PV expansion can be done up to the technical limits of the existing grid. 
\begin{equation}\label{eq:MaxProd}
    0 \leq p_{h,pv} \leq \sum_{gi(pv,i)} \overline{P}_{pv} \cdot CF_{h,i,pv} \cdot scal \quad \forall h,pv
\end{equation}
\begin{equation}\label{eq:Gen+Curt}
    p_{h,pv} + sp_{h,pv} = \sum_{gi(pv,i)} \overline{P}_{pv} \cdot CF_{h,i,pv} \cdot scal \quad \forall h,pv
\end{equation}
Constraints~\eqref{eq:DetBin}-\eqref{eq:Regulation3} limit grid-active feed-in by the maximum power and prevent curtailment beyond this. Constraint~\eqref{eq:DetBin} sets the binary variable $\alpha_{h,i}$ to 1 as soon as PV generation ($\overline{P}_{pv} \cdot PF_{h,i,pv} \cdot scal$) minus consumption ($D_{h,i}^P$) exceeds the permitted feed-in power ($FL \cdot \overline{P}_{pv} \cdot scal$). If the binary variable is set, \eqref{eq:Regulation1} and \eqref{eq:Regulation2} guarantee that the PV production $p_{h,pv}$ corresponds exactly to the maximum permissible feed-in power ($FL \cdot \overline{P}_{pv} \cdot scal + D_{h,i}^P$) and~\eqref{eq:Regulation3} is not binding. Vice versa, when no curtailment is necessary,~\eqref{eq:Regulation3} sets the curtailed power $sp_{h,pv}$ to zero and thus guarantees that all the available energy is also generated. In the equations provided, $M$ represents a large positive constant. Its purpose is to enable logical implications in mixed-integer programming (MIP) models by activating or deactivating constraints based on the value of binary decision variables ($\alpha_{h,i}$).

\begin{equation}\label{eq:DetBin}
    \sum_{gi(pv,i)} \Big(\overline{P}_{pv} \cdot CF_{h,i,pv} \cdot scal \Big) - D_{h,i}^P \leq  FL \cdot \overline{P}_{pv} \cdot scal - eps + (M +eps) \cdot \alpha_{h,i} \quad \forall h,i
\end{equation}
\begin{equation}\label{eq:Regulation1}
    \sum_{gi(pv,i)} \Big(p_{h,pv} - FL \cdot \overline{P}_{pv} \cdot scal\Big) - D_{h,i}^P \leq + M \cdot (1 - \alpha_{h,i}) \quad \forall h,i
\end{equation}
\begin{equation}\label{eq:Regulation2}
    \sum_{gi(pv,i)} \Big(p_{h,pv} - FL \cdot \overline{P}_{pv} \cdot scal\Big) - D_{h,i}^P \geq - M \cdot (1 - \alpha_{h,i}) \quad \forall h,i
\end{equation}
\begin{equation}\label{eq:Regulation3}
    \sum_{gi(pv,i)} sp_{h,pv} \leq M \cdot \alpha_{h,i} \quad \forall h,i
\end{equation}

\begin{table}[width=.55\linewidth,cols=4,pos=h]
\caption{List of parameters used in the algorithm}
\label{tbl:Parameter}
\begin{tabular*}{\tblwidth}{@{} LLLR@{} }
\toprule
Parameter                &               & Unit              & Value     \\
\midrule
Costs of imports         & $C^{Imp}$     & €/MWh, €/MVArh    & 200       \\
Costs of exports         & $C^{Exp}$     & €/MWh, €/MVArh    & 200       \\
Costs of PNS/QNS         & $C^{PNS/QNS}$ & €/MWh, €/MVArh    & 100\,000  \\
Costs of EPS/EQS         & $C^{EPS/EQS}$ & €/MWh, €/MVArh    & 200\,000  \\
\bottomrule
\end{tabular*}
\end{table}

\section{Case Study Description}
\label{sec:Case-Study}
In the following section, the individual case studies are presented in detail. Section~\ref{subsec:Data-and-general-assumptions} explains general assumptions and describes the data. Section~\ref{subsec:descrip-networks} analyzes with the medium- and low-voltage grids.

\subsection{Data and general assumptions}
\label{subsec:Data-and-general-assumptions}

For the case study in this paper, real grid data/configurations, including associated generators and consumers, were used. Data was provided by the local distribution system operator of the Federal State of Lower Austria\footnote{Netz Niederösterreich GmbH}. Analyses for the evaluation of dynamic feed-in limitation were carried out for both medium- and low-voltage grids using the second-order cone approximation of the AC-OPF implemented LEGO.

For the medium-voltage cases, three representative medium-voltage grids were selected by the DSO, each with different characteristics: a rural grid, an urban grid, and a hybrid grid with both urban and rural characteristics.
For the low-voltage cases, a grid with a high penetration rate of smart meters and, therefore, a high-quality database was selected. All data was anonymized. A detailed description of the example grids is provided in Section~\ref{subsec:descrip-networks}.

The central objectives of the calculations for the three \textbf{medium-voltage} grids include the systematic identification of bottlenecks and the precise determination of the maximum installable PV capacity (in addition to the existing systems). In this context, only the worst-case scenario (1 hour) was considered, which assumes low consumption with simultaneous maximum generation. In an initial step, an unlimited expansion of additional PV systems up to the technical capacity limits of the grid was carried out. In coordination with the DSO, a distinction was made between scalable and non-scalable systems, for example, full feed-in systems that have already reached the maximum expansion capacity. In addition, the possibility of adding PV at nodes where no PV system currently exists was included. In order to avoid disadvantaging PV systems due to unfavorable connection points (such as a long distance to the next transfer point to a higher grid level or grid area with a high voltage increase), the uniform scaling factor, $scal$, was optimized for all scalable PV systems. This factor invests equally in all PV systems in the entire grid until expansion is stopped due to technical limits. For the case studies, the maximum feed-in power is reduced incrementally, i.e., 100\% of the installed PV power may be fed in, then 90\%--80\%--70\%, to investigate the specific effect on the expansion options and resulting grid bottlenecks in detail. To verify the precision of the model and the data utilized, an operating case was examined for each of the three medium-voltage grids, employing the DSO's power flow simulation program.

In addition to maximizing generation capacity, analyses of the \textbf{low-voltage} grid also focus on total generation and the potential for alleviating generation peaks. Here, a full year with hourly resolution was analyzed in the optimization model. In contrast to the medium-voltage cases, however, it is not possible to expand existing systems. Instead, the construction of new PV systems was only possible for households without existing systems.

In order to be able to cover multiple regulatory developments for medium and low voltage, a distinction was made between two cases: (a) only new PV installations may be curtailed by the dynamic feed-in limitation; and (b) in addition to the new PV installations, existing PV installations are also limited in their feed-in. For an actual implementation of a dynamic feed-in limitation in Austria, case~(a) would be more likely, as it would be very challenging to change existing  feed-in contracts.

\subsection{Description of the grids}
\label{subsec:descrip-networks}

Three characteristic medium-voltage grids were analyzed: a rural grid, an urban grid, and a hybrid grid with both rural and urban characteristics. In addition, a low-voltage grid was also examined.

To systematically analyze the impact of feed-in limitation and a comparison of the different levels of dynamic feed-in limitation, PV expansion is pushed to the technical limits of the grid. This concerns both the thermal limit of power lines/cables and voltage limits, which must be strictly adhered to at all nodes in the grid. In the medium-voltage grids, the lower voltage limit is 0.95~p.u. and the upper limit is 1.03~p.u. For expansion, the existing PV capacity is used as the basis for expansion at nodes to which a scalable system is connected. For nodes with only non-scalable systems and without any consumption, there is no possibility for further expansion, as these are already full feed-in systems at maximum expansion. At nodes with consumption, but without existing systems, the mean value of all scalable PV systems was used as scalable capacity at these nodes. These candidates are adapted to the technical limits of the medium-voltage grid by the scaling variable, which can take a continuous value between zero and infinity. In addition, a sensitivity analysis is carried out to investigate the influence of consumption on the expansion in the worst-case scenario. To this end, all scenarios in the medium-voltage range were calculated both with the specified consumption and with a 10\% or 20\% increase in consumption. Imports/exports to the higher-level grids take place via 110/20~kV transformers, respectively. These transformers are located upstream of Node 0 and are not explicitly shown in the subsequent Figures. 
In the first three rows of Table~\ref{tbl:Grid-Infos}, the main information of the three medium-voltage grids is shown.
There are a total of 260 metering points in the low-voltage grid.  For simplicity, metering points in multi-party buildings have been combined into one node. In order to add new PV systems, a scalable PV system was added at nodes with consumption but no existing PV system. The generation profile for the new systems is based on full feed-in in close proximity to the low-voltage grid.
Unlike the medium-voltage grids, the low-voltage grid had a complete hourly time series (8760 hours) for each consumer, and the existing generation units are based on actual metering data provided by the DSO. A voltage range of 0.9~p.u. and 1.1~p.u. is specified as a limiting factor for the expansion of new PV systems, in addition to the thermal limiting power of the lines. All other key factors are listed in the last row of Table~\ref{tbl:Grid-Infos}.
It should be noted that the transformers in medium- and low-voltage grids are not to be regarded as bottlenecks as they can be replaced by more powerful ones if required.

\begin{table}[width=1\linewidth,cols=4,pos=h]
\caption{Key information of the analyzed grids}
\label{tbl:Grid-Infos}
\begin{tabular*}{\tblwidth}{@{} L|R|R|R|R|R|RR|RR|RR|RR|RR@{} }
\toprule
&\makecell{Trans-\\former} & Voltage &\makecell{Total\\Length} &Nodes &\makecell{Total Energy\\Demand} &\multicolumn{2}{L|}{RoR} &\multicolumn{2}{L|}{Wind} &\multicolumn{2}{L|}{PV$_{scal}$} &\multicolumn{2}{L|}{PV$_{non\:scal}$} &\multicolumn{2}{L}{Fossil}\\
          &  MVA &  kV & km  & \#    &   MWh    & \# & kW  & \# & kW  & \# & kWp   & \# & kWp   & \# & kW  \\
\midrule
Rural MV  &   25 &  20 & 21  & 158   &   1.20  & 1  & 4.4 & 3  & 320 & 53 & 6\,838 &  4 & 1\,432 & -  & -   \\
Urban MV  &   40 &  20 & 4   & 110   &   2.00  & 2  & 208 & -  & -   & 27 & 2\,028 &  2 & 2\,622 & 1  & 300 \\
Hybrid MV &   40 &  20 & 78  & 267   &   2.51  & 1  & 80  & -  & -   & 74 & 8\,264 &  1 &    300 & -  & -   \\
LV        & 0.25 & 0.4 & 6.9 & 180   & 483.72  & -  & -   & -  & -   &  - & -      & 39 &    352 & -  & -   \\
\bottomrule
\end{tabular*}
\end{table}

\section{Results and Discussion}
\label{Sec:Results}
In this section, we present the results for the medium- and low-voltage grids described in section~\ref{sec:Case-Study}. Section~\ref{subsec:Results-medium-voltage-grids} contains the results of the medium-voltage grids. The results of the low-voltage grid can be found in section~\ref{subsec:Results-low-voltage-grid}.

\subsection{Medium-voltage grids}
\label{subsec:Results-medium-voltage-grids}
In the following sections~\ref{subsubsec:Results-rural-medium-voltage-grid}-\ref{subsubsec:Results-hybrid-medium-voltage-grid}, the results of the three analyzed medium-voltage grids and the limitations for the integration of additional PV systems are discussed in detail.

\subsubsection{Rural medium-voltage grid}
\label{subsubsec:Results-rural-medium-voltage-grid}
\begin{figure}[h!]
\centering
\includegraphics[width=1\textwidth]{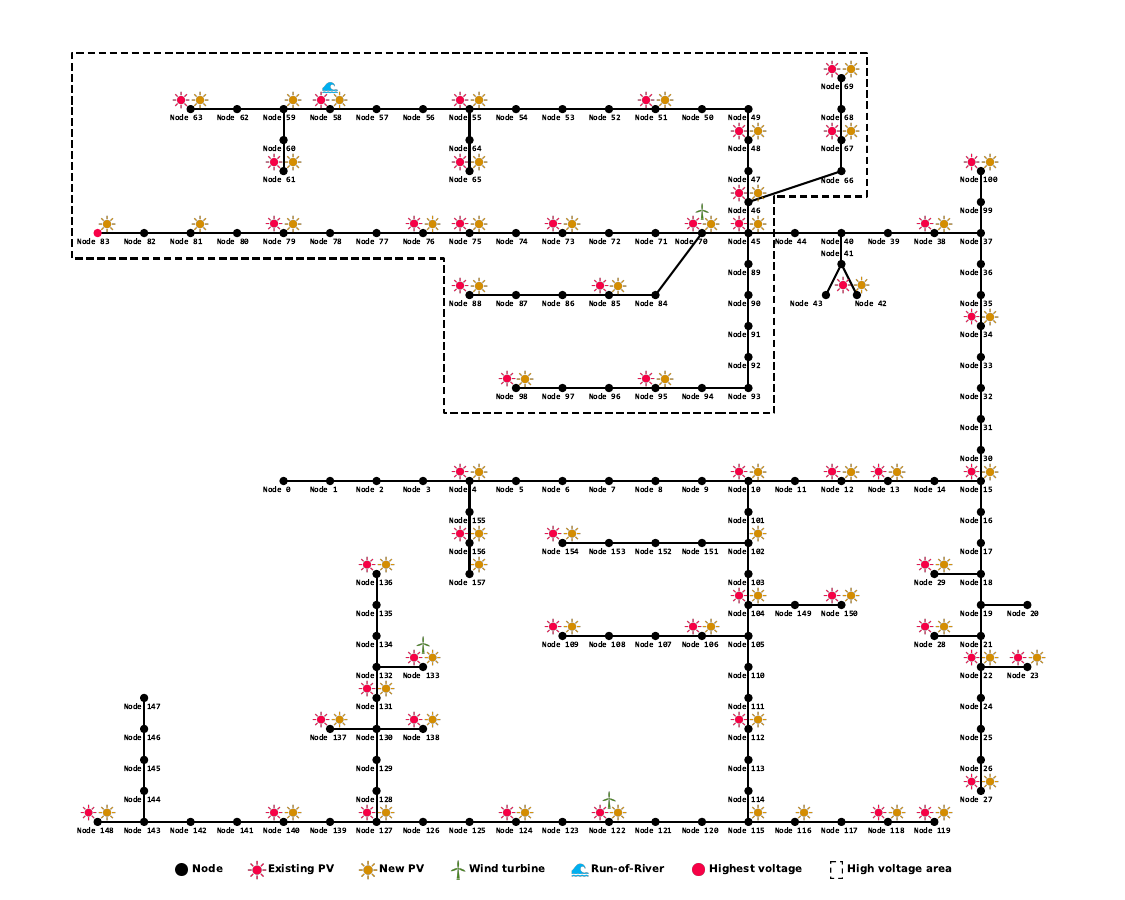}
\caption{Schematic representation of the rural medium-voltage grid with existing and new generators, as well as the bottlenecks for further expansion.}\label{fig:Rural-MV_Grid_Res}
\end{figure}

The rural medium-voltage grid under consideration has a total length of 21~km and contains 61 generators, including three wind turbines (green wind turbine symbol) and one run-of-river power plant (blue wave symbol). The remaining 57 generators are existing PV systems (red sun symbol), of which 53 are expandable (yellow sun symbol). Additional details regarding the grid can be found in Table~\ref{tbl:Grid-Infos} and a schematic illustration of the radial grid is provided in Figure~\ref{fig:Rural-MV_Grid_Res}. To evaluate the potential for additional PV integration under different curtailment strategies, two cases are examined:

Case~(a): Assuming that only newly built PV systems can be curtailed according to the modeled PV feed-in limitation, in this rural medium-voltage grid, between 1.05~MWp and 1.07~MWp of PV capacity could be added, depending on the feed-in limitation (between 70\% and 100\% of the installed module peak power). This is depicted in Figure~\ref{fig:Res_Rural-MV-Power} (dark red bars). Hence, hardly any improvement with respect to additional PV integration could be achieved despite the dynamic feed-in limitation. This is due to the existing PV systems in the area of nodes 45 to 98, which already cause increased voltage in this area. As these existing systems cannot be curtailed in this analysis, the dynamic feed-in limitation could not have any effect due to the consumption at the corresponding nodes. 

Case~(b): For the consideration of possible curtailment of all PV systems (old and new), the sum of the maximum installable module peak power of new PV systems in the entire system could be increased from 1.05~MW to 3.7~MWp through dynamic feed-in limitation. The dark green bars in Figure~\ref{fig:Res_Rural-MV-Power} show the exact effect of the feed-in limitation on installable module peak power. The limiting factor for this analysis was also the voltage in the range of nodes 45 to 98 and reaching the upper voltage limit of 1.03~p.u. at node 83. 

Figure~\ref{fig:Rural-MV_Grid_Res} illustrates the schematic representation of the rural medium-voltage grid under consideration, with node 83 (highlighted red) marked as a critical point in the grid that has reached the maximum permissible voltage of 1.03~p.u. The black dashed area covers all other nodes that also fall within the upper range of the permitted voltage band ($\geq$1.028~p.u.).

In order to investigate the influence of consumption on the installable module peak power, the demand (1.2~MW) was increased by 10\% and 20\% for both cases (a) and (b). Figure~\ref{fig:Res_Rural-MV-Power} shows the influence on output in light red if only new systems are regulated and in light green if both new and existing systems can be curtailed. The graph on the right shows that the higher consumption has a progressively smaller effect as the installed peak module power in the grid increases. For example, an increase of 22\% could be achieved for 100\% possible feed-in, whereas the effect is limited to just under 2.8\% for a curtailment to 70\%.

The expansion of PVs in this rural medium-voltage grid is limited by the voltage limit. The additional capacities that can be achieved through dynamic feed-in limitation are low due to the existing systems and the associated high voltage increase. For example, the installed peak module power for case~(a) could only be increased by 2\%. Additional regulation of the existing systems, on the other hand, increased the expansion by around 250\%. A higher consumption in the grid also increases the maximum module peak power that can be installed. However, this higher consumption has a progressively smaller impact as the effectiveness of dynamic feed-in limitation increases.

\subsubsection{Urban medium-voltage grid}
\label{subsubsec:Results-urban-medium-voltage-grid}
\begin{figure}[h!]
\centering
\includegraphics[width=\textwidth]{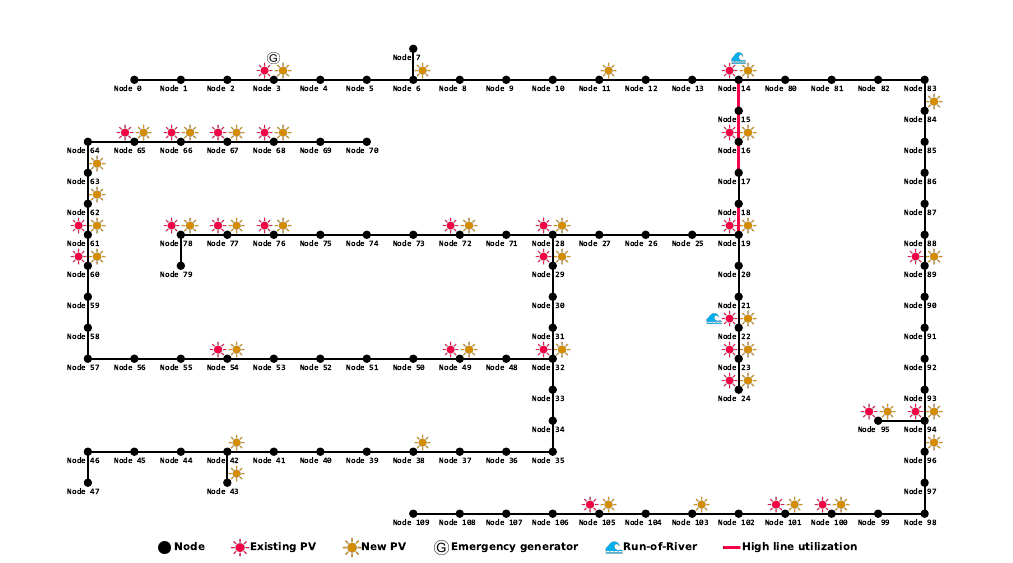}
\caption{{Schematic representation of the urban medium-voltage grid with existing and new generators, as well as the bottlenecks for further expansion.}}\label{fig:Urban-MV_Grid_Res}
\end{figure}

The urban medium-voltage grid under consideration has a total length of only 4~km and contains 32 generators. Of these, two are run-of-river power plants (blue wave symbol) and one is an emergency power plant (black generator symbol) belonging to a company that is routinely tested on the grid. The remaining 29 generators are existing PV systems (red sun symbol), 27 of which are expandable (yellow sun symbol). Further details regarding the grid can be found in Table~\ref{tbl:Grid-Infos} and a schematic illustration of the radial grid is given in Figure~\ref{fig:Urban-MV_Grid_Res}. To evaluate the potential for additional PV integration under different curtailment strategies, two cases are examined:

Case~(a): For the urban medium-voltage grid, assuming that only new PV systems are limited to a grid-active feed-in power of between 70\% and 100\%, between 14.8~MWp and 18.5~MWp of module peak power could be added. The dark red bars in Figure~\ref{fig:Res_Urban-MV-Power} show the effect of dynamic feed-in limitation on the installable capacity. In all cases, the thermal power limits of the line sections between nodes 14 and 17 as well as between nodes 18 and 19 were the limiting factors.

Case~(b): For the possibility of curtailment of all PV systems, the maximum installable module peak power could be increased from 14.8~MW to 21~MWp using dynamic feed-in limitation. Figure~\ref{fig:Res_Urban-MV-Power} (dark green bars) show the effect of feed-in limitation on the installable module peak power as a function of the permissible curtailment. Analogously to case~(a), here, the line sections between nodes 14 and 17 and the section between nodes 18 and 19 limit further expansion.

The thermal power limits of some line sections (power bottleneck) between nodes 14 and 17, as well as between nodes 18 and 19, emerged as a limiting factor for the urban medium-voltage grid. All line sections between nodes 14 and 19 have a thermal power limit of 11~MW, the only exception being the section between nodes 17 and 18, which has already been replaced by more powerful 15.2~MW cabling. Figure~\ref{fig:Urban-MV_Grid_Res} illustrates the identified bottlenecks, marked red in the schematic diagram.

Figure~\ref{fig:Res_Urban-MV-Power} shows the effects of a 10\% and 20\% increase in consumption. The light red {bars} show the installable peak module outputs for case~(a) and the light green bars, those for case~(b). The difference due to the increase in consumption between the three graphs is marginal because of the low load of 2.5~MW compared with the installed capacity of at least 19~MW (existing systems 4.33~MW plus new systems). The effect of higher consumption is also reduced in the urban medium-voltage grid as the installed peak module power increases.

\subsubsection{Hybrid medium-voltage grid}
\label{subsubsec:Results-hybrid-medium-voltage-grid}
\begin{figure}[h!]
\centering
\includegraphics[width=1\textwidth]{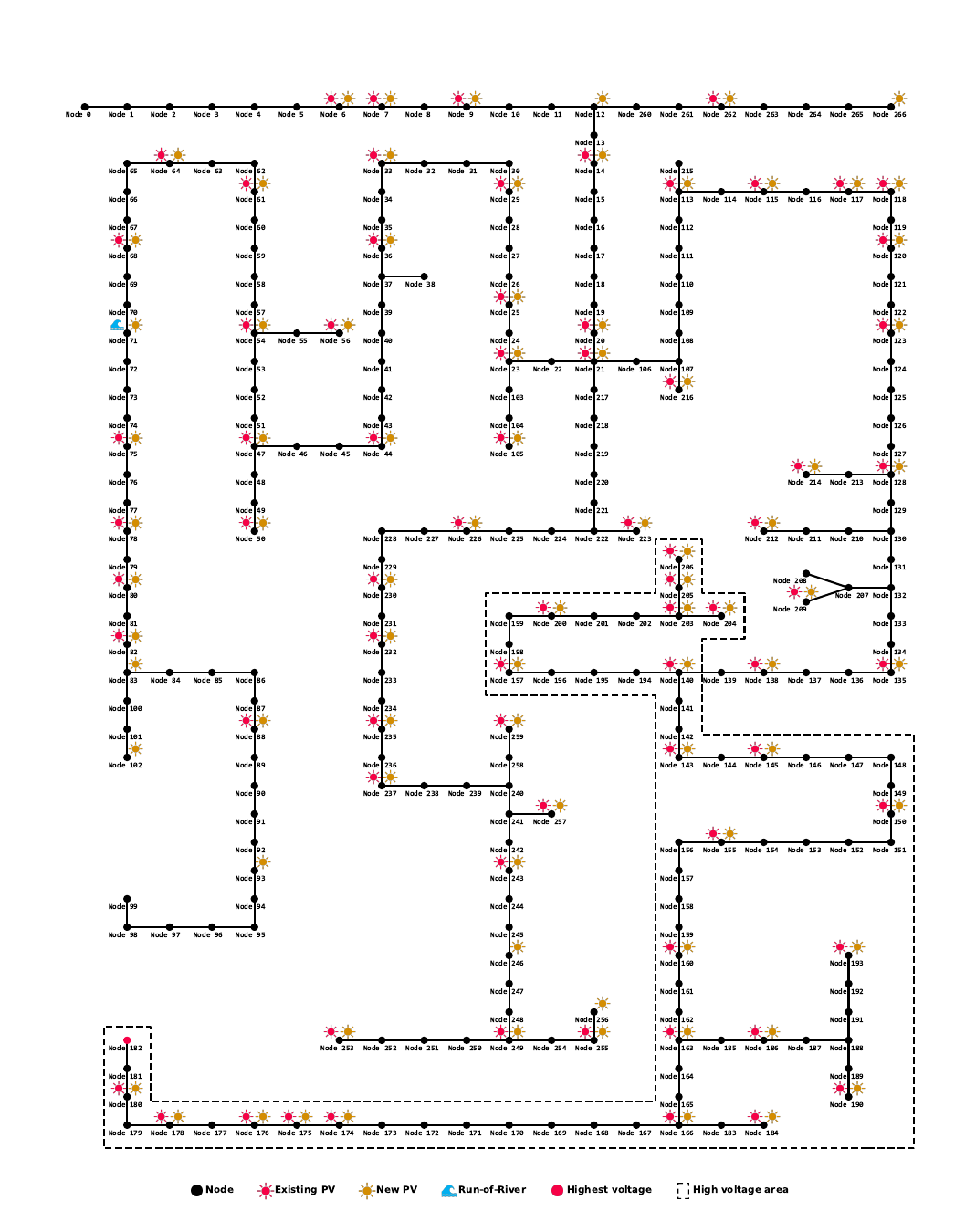}
\caption{Schematic representation of the hybrid medium-voltage grid with existing and new generators, as well as the bottlenecks for further expansion.}\label{fig:Hybrid-MV_Grid_Res}
\end{figure}

\begin{figure}
\centering
\begin{minipage}{.32\textwidth}
  \centering
  \includegraphics[width=1\linewidth]{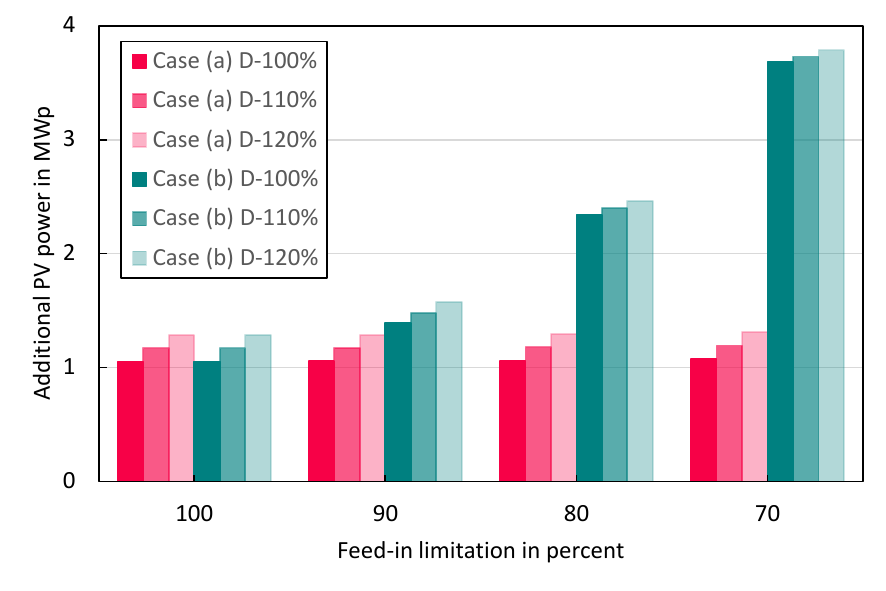}
  \subcaption{Rural Grid}
  \label{fig:Res_Rural-MV-Power}
\end{minipage}
\begin{minipage}{.32\textwidth}
  \centering
  \includegraphics[width=1\linewidth]{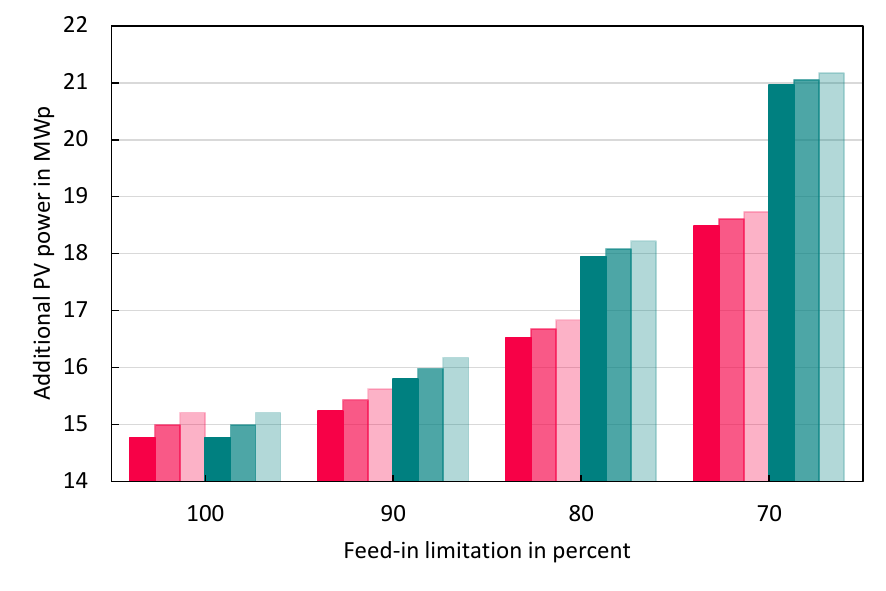}
  \subcaption{Urban Grid}
  \label{fig:Res_Urban-MV-Power}
\end{minipage}
\begin{minipage}{.32\textwidth}
  \centering
  \includegraphics[width=1\linewidth]{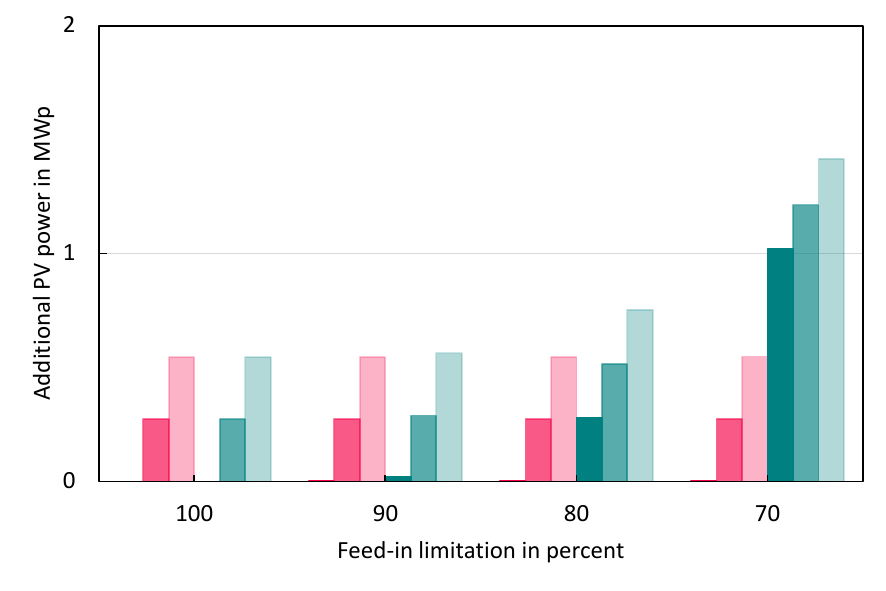}
  \subcaption{Hybrid Grid}
  \label{fig:Res_Hybrid-MV-Power}
\end{minipage}
\caption{Maximum installable power of new PV systems depending on the maximum feed-in power in medium-voltage grids.}
\label{fig:Res_MV-Power}
\end{figure}

The hybrid medium-voltage grid, which exhibits characteristics of both rural and urban grids, has a total length of 78~km and contains 76 generators, including one run-of-river power plant (blue wave symbol) and 75 existing PV plants (red sun symbol). All except one of the PV plants can be expanded (yellow sun symbol). Further details regarding the grid can be found in Table~\ref{tbl:Grid-Infos}, and a schematic illustration of the radial grid is provided in Figure~\ref{fig:Hybrid-MV_Grid_Res}. To evaluate the potential for additional PV integration under different curtailment strategies, two cases are examined:

Case~(a): Assuming that only new PV systems can be curtailed, no additional systems could be integrated into the grid, as the voltage cap was already reached without any new additions. As can be seen in Figure~\ref{fig:Res_Hybrid-MV-Power} (dark red bars), no improvement could be achieved through dynamic power regulation either. This is due to the existing systems in the area of nodes 139 to 206, which already cause a higher voltage increase in this area. Since these existing systems could not be curtailed in this analysis, the dynamic feed-in limitation could not have any effect due to the consumption at the corresponding nodes.

Case~(b): For the consideration of possible curtailment of all PV systems, the maximum installable module peak power could be increased from 0~MWp to just over 1~MWp through dynamic feed-in limitation. The dark green bars in Figure~\ref{fig:Res_Hybrid-MV-Power} show the effect of the feed-in limitation on the installable module peak power. The limiting factor for this analysis was also the voltage in the area of nodes 139 to 206 and reaching the upper voltage limit of 1.03~p.u. at node 182.

Figure~\ref{fig:Hybrid-MV_Grid_Res} shows the schematic representation of the medium-voltage grid under consideration, with node 182 (highlighted red) marked as a critical point in the grid that has reached the maximum permissible voltage of 1.03~p.u. The black dashed area includes all other nodes that also have a high voltage increase.

\clearpage

The influence of consumption on the installable module peak power is shown in Figure~\ref{fig:Res_Hybrid-MV-Power}, the light red bars show the influence on power if only new systems are regulated, and the light green bars if both new and existing systems can be curtailed. For this purpose, consumption (2.5~MW) was increased by 10\% and 20\% for both cases (a) and (b). For case~(a), the increase in consumption can also increase the additional installable capacity, but the dynamic feed-in limitation has no effect due to the negligible PV expansion. Case~(b) shows that the increase in consumption with a decreasing percentage of feed-in power has a more minimal impact on the growth of supplementary PV systems.

The expansion of PVs in this hybrid medium-voltage grid, which is characterized by both rural and urban areas, is limited by the upper voltage limit of 1.03~p.u. The additional capacities that can be achieved through dynamic feed-in limitation are low due to the existing systems and the associated high voltage increase. The installed peak module capacity could not be increased for case~(a), but for case~(b) the expansion could be increased to at least 1~MWp. A higher consumption in the grid also increases the maximum module peak power that can be installed. However, the higher consumption has a progressively smaller impact as the effectiveness of dynamic feed-in limitation increases.

\subsection{Low-voltage grid}
\label{subsec:Results-low-voltage-grid}
The low-voltage grid under consideration extends over a total length of 6.9 km and includes 260 metering points as well as 39 existing PV systems (red sun symbol). The expansion of new PV systems (yellow sun symbol) is a viable option at all nodes where consumption is present but an existing PV system is not installed. Further details about the grid are provided in Table~\ref{tbl:Grid-Infos}, while a schematic illustration can be found in Figure~\ref{fig:LV-Grid-Res}. The ensuing sections present the results of the analysis. Subsection~\ref{subsubsec:LV-installable-capacity} examines the impact of dynamic feed-in limitation on the maximum installable PV capacity, while subsection~\ref{subsubsec:LV-energy-increase} focuses on the resulting energy generation. Subsection~\ref{subsubsec:LV-losses} analyzes curtailment losses as a function of the maximum feed-in, and finally, subsection~\ref{subsubsec:LV-bottlenecks} identifies bottlenecks within the grid.


\subsubsection{Installable capacity}
\label{subsubsec:LV-installable-capacity}
Case~(a): By using the dynamic feed-in limitation, a maximum grid-active feed-in power ranging between 70\% and 100\% of the installed peak module power is set. This methodology facilitates the augmentation of the installable module peak power for new systems. In the absence of feed-in limitation (100\%), the installable module peak power is 9.4 kWp. However, when a feed-in limitation of 70\% is implemented, the potential for increase in installable module peak power reaches 12.4 kWp. This represents a 32\% increase in peak installed PV capacity. Figure~\ref{fig:LV_Power-per-generator} (red bars) shows the correlation between the maximum installable module peak power per system and the percentage of permissible feed-in to the low-voltage grid.

Case~(b): With the additional curtailment of existing systems, the installable system size could be increased from 9.4~kWp without any feed-in limitation (100\%) up to 12.6~kWp for a limitation of 70\%, which corresponds to an increase of 34\%. The other values for a curtailment of 80\% or 90\% can be seen in Figure~\ref{fig:LV_Power-per-generator}. In both cases, the expansion was limited by the thermal limit of the power line between node 60 and the low-voltage side of the transformer (node~0).

\subsubsection{Increase in energy generation}
\label{subsubsec:LV-energy-increase}
Case~(a): As a complete time series was available for the entire low-voltage grid for a whole year (8760 hours), it was possible to quantify not only the capacity but also the yearly generated energy that was gained through dynamic feed-in limitation. Up to 1.28~GWh (generated by PV systems) could be generated due to the dynamic curtailment of new systems only -- without any limitation of feed-in capacity (100\%). In comparison, the energy generation was 1.02~GWh per year. The red bars in Figure~\ref{fig:LV_Generated-energy} show the correlation between the total generation of all systems and the grid-active feed-in power as a function of system output. The percentage increase in total reaches up to 25.4\% more than without dynamic feed-in limitation for this grid.

Case~(b): For the limitation of feed-in systems for all PV systems, generation could be increased to up to 1.3~GWh per year with a grid-active feed-in power of 70\% of module peak power. This corresponds to an increase of almost 27\% compared with the possibility of being able to feed everything into the grid. As illustrated in Figure~\ref{fig:LV_Generated-energy}, the green bars demonstrate the relationship between total generation and permissible feed-in power.

\subsubsection{Losses due to dynamic feed-in limitation}
\label{subsubsec:LV-losses}
\begin{figure}[h!]
\centering
\begin{minipage}{.32\textwidth}
  \centering
  \includegraphics[width=1\linewidth]{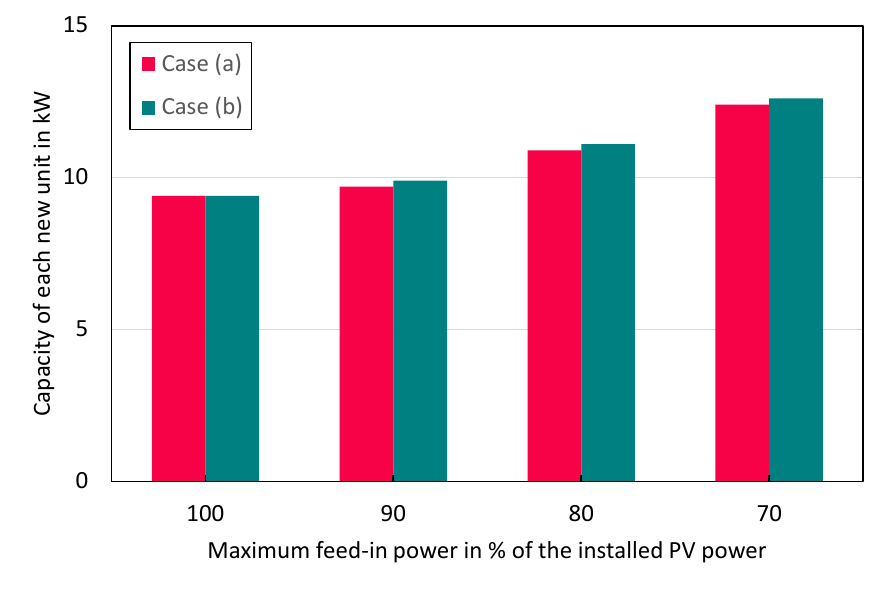}
  \subcaption{Maximum installable capacity}
  \label{fig:LV_Power-per-generator}
\end{minipage}
\begin{minipage}{.32\textwidth}
  \centering
  \includegraphics[width=1\linewidth]{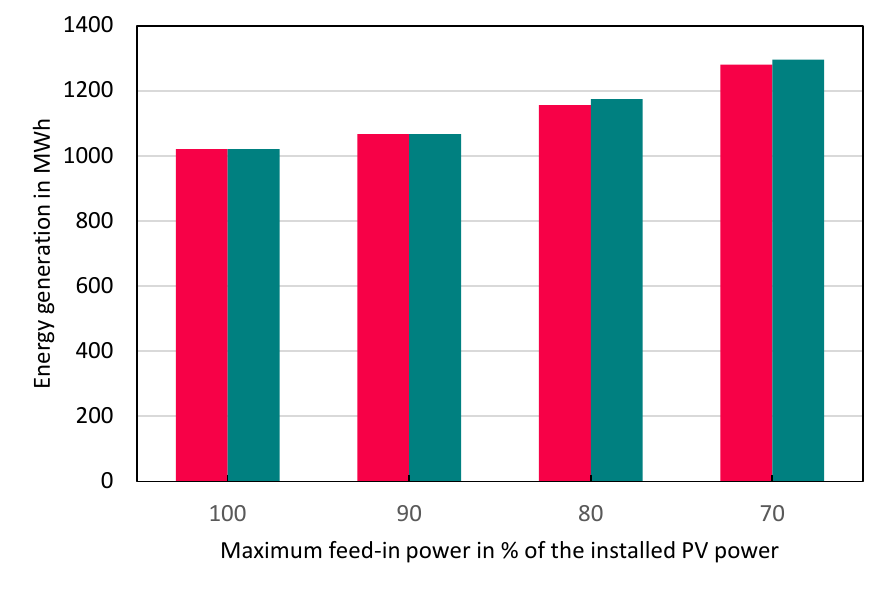}
  \subcaption{Total generation in the grid}
  \label{fig:LV_Generated-energy}
\end{minipage}
\begin{minipage}{.32\textwidth}
  \centering
  \includegraphics[width=1\linewidth]{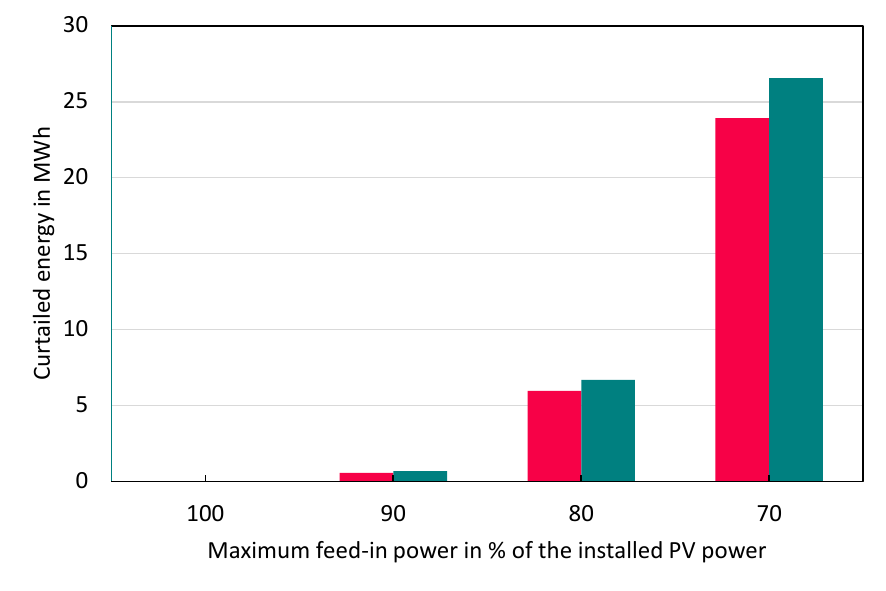}
  \subcaption{Total curtailment in the grid}
  \label{fig:LV_Curtailed-energy}
\end{minipage}
\caption{Maximum installable power of new PV systems, total energy generation, and curtailed energy} depending on the maximum feed-in power in low-voltage grids.
\label{fig:Res_LV}
\end{figure}

Case~(a): Limiting the grid-active feed-in power of the system output also results in the curtailment of energy from PV systems. As can be seen in Figure~\ref{fig:LV_Curtailed-energy} (red bars), up to 24~MWh per year were lost in the entire low-voltage grid under consideration for the curtailment of new systems only. This corresponds to slightly more than 1.8\% of the maximum possible annual energy generation.

Case~(b): If curtailment is allowed in all systems, in comparison with case~(a), with a maximum feed-in percentage of 70\%, the curtailed energy increases by over 2.5~MWh to just under 26.6~MWh. This corresponds to 2\% of the maximum possible annual energy generation. Figure~\ref{fig:LV_Curtailed-energy} (green bars) also shows the values for 80\% and 90\% feed-in limitation.\\

\begin{figure}[h!]
\centering
\includegraphics[width=84mm]{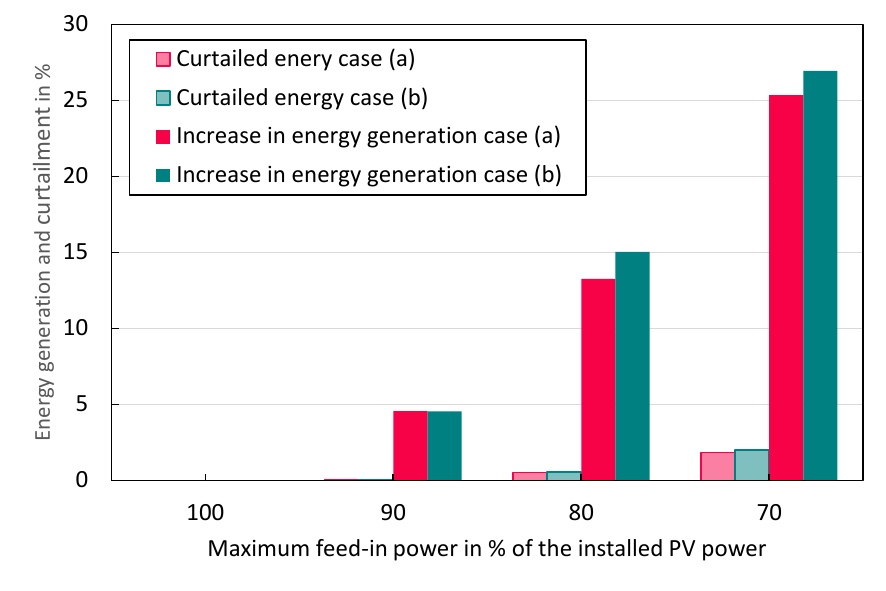}
\caption{Comparison of generated and curtailed energy in the low-voltage grid}\label{fig:LV-increase-losses}
\end{figure}

As illustrated in Figure~\ref{fig:LV-increase-losses}, for both cases (a) and (b), there is a direct comparison between the additional generation and the curtailed energy.

\subsubsection{Bottleneck in the low-voltage grid}
\label{subsubsec:LV-bottlenecks}
\begin{figure}[h!]
\centering
\includegraphics[width=\textwidth]{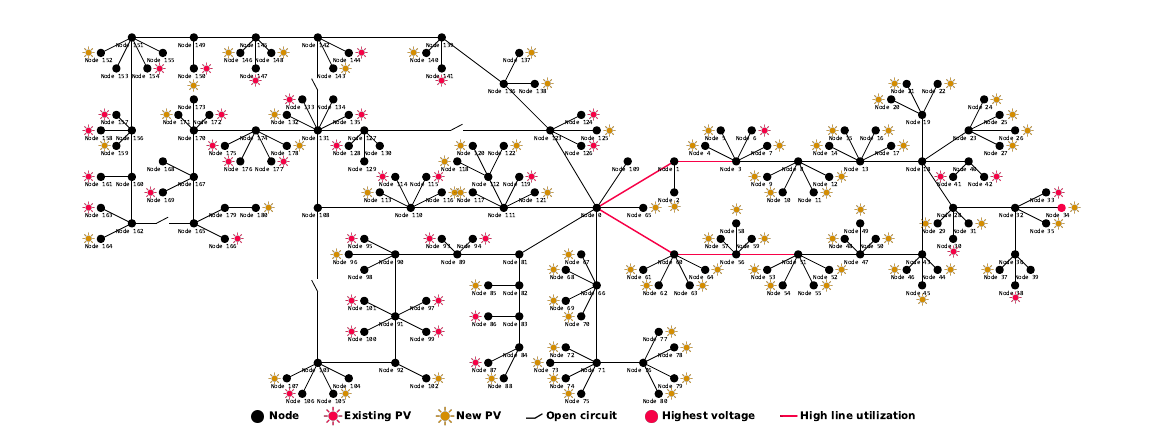}
\caption{Schematic representation of the low-voltage grid with existing and new PV systems, as well as the bottlenecks for further expansion.}\label{fig:LV-Grid-Res}
\end{figure}

The limiting factor (bottleneck: power) for further expansion in all low-voltage model runs was the maximum power of the line between node 60 and node 0 (low-voltage side transformer). However, the connections upstream (node 60 - node 56 - node 51 - node 47) and on the other side of the conductor loop (node 0 - node 1 - node 3) also exhibited high utilization. The maximum voltage occurred at node 34, but at just +6\% of the nominal voltage the value remained well below the voltage limits of ±10\%. Figure~\ref{fig:LV-Grid-Res} shows the lines with the highest utilization and the node with the highest voltage increase.\\

\clearpage
\section{Conclusions}
\label{sec:Conclusion}

In this paper, we proposed a mathematical formulation to analyze the effects of dynamic feed-in limitations of PV systems on distribution grids, which is currently being discussed in Austria. This method improves the potential installed capacity of PV systems on existing grids, thereby enabling a more expeditious expansion of renewable energies. We validated the model with real examples of Austrian medium- and low-voltage grids. By limiting grid-active feed-in power, additional customers could be granted access to the grid more quickly compared with having to wait for grid expansion. We analyzed three medium-voltage grids with different characteristics (rural, urban, hybrid) and one low-voltage grid using the LEGO expansion planning model with dynamic feed-in limitation extension. With this model extension, we have determined the maximum possible expansion of new PV systems.

The results from the two medium-voltage grids (rural and hybrid) that already have a relevant installed capacity of PV indicate that only the curtailment of new systems does not allow additional expansion of new PV systems. In comparison, the curtailment of both new and existing systems enables a significantly higher expansion. In the urban medium-voltage grid and the low-voltage grid, the curtailment of only new systems has demonstrated that this also has a significant (up to 32\%) influence on expansion, due to the as yet moderate number of existing PV systems. The curtailment of all systems has also further increased the potential PV expansion in these grids. 
Irrespective of the limitation of expansion (reaching the upper voltage or maximum thermal power limit), it is shown that in grids with an already high penetration rate of PV systems, only the curtailment of all systems brings a noticeable improvement. In grids with sufficient available grid capacities, both cases of curtailment have a positive effect on the potential capacity.
The annual analysis of the low-voltage grid shows that the minor loss (2\%) of generation caused by limiting feed-in power is relatively small in comparison with the benefits of the resulting freed-up grid capacities and the associated additional PV system generation.

A limitation of the model pertains to the distribution of new PV systems, which poses a significant challenge to forecasting. Additionally, the model does not consider variations in the orientation of PV systems or the effects of shading, both of which can have a substantial impact on generation capacity. Another limitation is the time resolution of one hour, which is based on average values for demand and PV generation. Achieving a higher temporal resolution would facilitate a more precise capture of demand and generation peaks, which are of significant importance for determining the maximum installable PV capacity.

Future research could include the  exploration of even more dynamic limitations of PV systems (e.g., accounting for network constraints day-ahead or in real-time, centralized or decentralized smart storage management) and the associated effects on the individual grid levels and PV integration, as well as what needs to be done to motivate PV system owners to participate.

In order to maintain the progressive expansion of renewable energies, particularly PV, and avoid a slowdown caused by inadequate or delayed grid expansion, the implementation of intelligent approaches is imperative to further accelerate the expansion process.



\printcredits

\bibliographystyle{elsarticle-num}




\end{document}